\documentclass[12pt,a4paper]{article}
\usepackage{amsfonts,amsmath,amssymb,enumerate,amsthm,ifthen,graphpap,longtable,time,xspace,graphicx}
\newcommand{\FigureSwitch}{on}
\newcommand{\CommentSwitch}{on}
\newcommand{\eop}{\hspace*{\fill}~$\square$} 
\theoremstyle{plain}
\numberwithin{equation}{section}
\newtheorem{theorem}{Theorem}[section]
\newtheorem{proposition}[theorem]{Proposition}
\newtheorem{lemma}[theorem]{Lemma}

\theoremstyle{definition}
\newtheorem{Remark}[theorem]{Remark}

\newcommand{\RootPath}{.}

\newcommand{\ExternalFiguresPath}{\RootPath/figures}

\newcommand{\natur}{\ensuremath{\mathbb{N}}}
\newcommand{\real}{\ensuremath{\mathbb{R}}}

\newcommand{\setcond}[2]{\left\{ #1 : #2 \right\}} 
\newcommand{\setcondbegin}[2]{\left\{ #1 : #2 \right.} 
\newcommand{\setcondend}[1]{\left.  #1 \right\}} 

\newcommand{\IncludeGraph}[2]{
\ifthenelse{\equal{\FigureSwitch}{on}}{
	\includegraphics[#1]{\ExternalFiguresPath/{#2}}
	}{
		\fbox{\texttt{#2}}
	}
}

\newcommand{\card}{\mathop{\mathrm{card}}\nolimits}

\newcommand{\intr}{\mathop{\mathrm{int}}\nolimits}

\newcommand{\MxL}{\left[} 
\newcommand{\MxR}{\right]} 

\newcommand{\ThmTitle}[2][]{\ifthenelse{\equal{#1}{}}{\emph{(#2)}}{\emph{(#2; #1)}}}
\newcommand{\notion}[2][]{\emph{#2}\xspace} 

\newcommand{\floor}[1]{\left\lfloor #1 \right\rfloor}
\newcommand{\ceil}[1]{\left\lceil #1 \right\rceil}

\newcommand{\eps}{\varepsilon}
\renewcommand{\card}{\#}

\newcommand{\modulo}[1]{ \ (\mathrm{mod} \, #1)}

\newcommand{\mycite}[2]{\ifthenelse{\equal{#2}{}}{\cite{#1}}{\cite[#2]{#1}}\xspace}

\newcommand{\RAGbook}[1][]{\mycite{MR1659509}{#1}}

\newcommand{\comment}[1]{\ifthenelse{\equal{\CommentSwitch}{on}}{/* #1 */}{}}

\newcommand{\edges}{E}
\newcommand{\rep}[3]{\{ #1 #3 0, \ldots, #2 #3 0 \} }

\begin{document}
\title{Description of Polygonal Regions by Polynomials\\ of Bounded Degree}
\date{\small \today}
\author{\small Gennadiy Averkov\footnote{Work supported by the German
Research Foundation within the Research Unit 468 ``Methods from
Discrete Mathematics for the Synthesis and Control of Chemical
Processes''.}\,\,, Christian Bey}
\maketitle

\begin{abstract} 
We show that every (possibly unbounded) convex polygon $P$ in $\real^2$ with $m$ edges can be represented by  inequalities
$p_1 \ge 0,\ldots,p_n \ge 0,$ where the $p_i$'s are products of at most $k$ affine functions each
vanishing on an edge of $P$ and $n=n(m,k)$ satisfies $s(m,k) \le n(m,k) \le (1+\eps_m) s(m,k)$ with
$s(m,k):=\max \{m/k,\log_2 m\}$ and $\eps_m \rightarrow 0$ as  $m \rightarrow \infty$.
This choice of $n$ is asymptotically best possible. An analogous result on representing the interior of $P$ in the form
$p_1 > 0,\ldots, p_n > 0$ is also given. For $k \le m/\log_2 m$ these statements remain valid for representations with arbitrary polynomials of degree not exceeding  $k$. 
\end{abstract} 

\newtheoremstyle{itsemicolon}{}{}{\mdseries\rmfamily}{}{\itshape}{:}{ }{}
\newtheoremstyle{itdot}{}{}{\mdseries\rmfamily}{}{\itshape}{.}{ }{}
\theoremstyle{itdot}
\newtheorem*{msc*}{2000 Mathematics Subject Classification} 

\begin{msc*}
  Primary: 14P10, 52B11;  Secondary: 52A10
\end{msc*}

\newtheorem*{keywords*}{Key words and phrases}

\begin{keywords*}
Gray code; polygon; polynomial; semi-algebraic set; Theorem of Br\"{o}cker and Scheiderer.
\end{keywords*}

\section{Introduction}

A set $P$ in $\real^d$ is said to be a \notion{$d$-dimensional polyhedron} if $P$ is intersection of finitely many closed half-spaces and has non-empty interior. Thus, a natural description of $P$ is given in terms of its $\mathcal{H}$-{\em representation}, i.e. a collection of non-strict affine inequalities. As was suggested in \cite{MR1976602}, \cite{MR2166533} and \cite{Henk06PolRep}, one could try to use more general polynomial representations of $P$, called $\mathcal{P}$-{\em representations}, for possibly developing new effective methods for solving combinatorial optimization problems. That is, one may look for representations of the form 
\begin{equation*}
	P = \{p_1 \ge 0, \ldots, p_n \ge 0 \} := \setcond{x \in \real^d}{p_1(x) \ge 0,\ldots,p_n(x) \ge 0},
\end{equation*} 
where $p_1,\ldots,p_n$ are real polynomials in $d$ variables. Naturally, one tries to find $\mathcal{P}$-representations of $P$ with a small number $n$ of polynomials.
Possible choices of $n$, also for the more general case when $P$ is a semi-algebraic set, can be derived using  results from real algebraic geometry, see \RAGbook, \cite{MR1393194}. It has been recently established that every $d$-dimensional polyhedron has a $\mathcal{P}$-representation with $d$ polynomials, see the end of this section for
references and further details.

In this note we are interested in small $\mathcal{P}$-representations of polyhedra consisting of polynomials of bounded degree.
We present a result in the $2$-dimensional case, namely, we determine asymptotically the minimal number of polynomials of degree at most $k$ which are
needed to describe a $2$-dimensional polyhedron with $m$ edges, for all degrees $k\le m/\log_2 m$. This improves earlier work by \cite{HenkMatzke}, see below. 

Let us call a $2$-dimensional polyhedron $P$ a \notion{polygon} (thus, we also allow unbounded polygons). Besides polynomial representations of $P$ we shall also investigate polynomial representations of the interior $\intr P$ of $P$ in terms of strict polynomial inequalities. Our main result is the following:

\begin{theorem} \label{main:cor} 
Let $m, k \in \natur$. Let $N(m,k)$ (resp. $\overline{N}(m,k)$) denote the minimal $n \in \natur$ such that for every polygon $P$
with $m$ edges there exist real polynomials $p_1,\ldots,p_n$ satisfying the conditions conditions (A) and (B) (resp.
($\overline{A}$) and (B)), where 
	\begin{enumerate} 
		\item[(A)] \label{repr:cond:open} $\intr P=\rep{p_1}{p_n}{>}$;
	        \item[($\overline{A}$)] \label{repr:cond:closed} $P=\rep{p_1}{p_n}{\ge}$;
		\item[(B)] \label{prod:cond:1:open} every $p_i$ has degree at most $k.$ 
	\end{enumerate}
Then for all $k$ satisfying $$k \le m / \log_2 m$$ one has  
	$$m/k \le N(m,k)\le  (1+ \eps_m) m/k,$$ where $\eps_m \rightarrow 0$ as $m \rightarrow +\infty$.\\
The same statement holds with $\overline{N}(m,k)$ instead of $N(m,k)$.
\end{theorem} 

The upper bound in Theorem \ref{main:cor} will be established by using polynomials $p_i$ each of which being a product of affine
functions. In fact, for such polynomials we are able to determine their minimal number $n$ asymptotically without
any restriction on $k$:

\begin{theorem} \label{main:thm}
Let $m, k \in \natur$. Let $N(m,k)$ (resp. $\overline{N}(m,k)$) denote the minimal $n \in \natur$ such that for every polygon $P$ with $m$ edges there exist real
polynomials $p_1,\ldots,p_n$ satisfying the conditions conditions (A) and (C) (resp.
($\overline{A}$) and (C)), where 
	\begin{enumerate} 
		\item[(A)] \label{repr:cond:open} $\intr P=\rep{p_1}{p_n}{>}$;
	        \item[($\overline{A}$)] \label{repr:cond:closed} $P=\rep{p_1}{p_n}{\ge}$;
		\item[(C)] \label{prod:cond:1:open} every $p_i$ is a product of at most $k$ affine functions each non-negative on
$P$ and vanishing on some edge of $P$. 
	\end{enumerate}
	Then, with $s(m,k) := \max \{m/k, \log_2 m \}$,
 $$s(m,k) \le N(m,k)\le \overline{N}(m,k) \le  (1+ \eps_m) s(m,k),$$
where $\eps_m \rightarrow 0$ as $m \rightarrow +\infty.$
\end{theorem} 

Prior to Theorem \ref{main:thm}, an upper bound $\overline{N}(m,k)\le m/k+\log_2 k$ was observed in \cite{HenkMatzke}
(see also \cite[Satz\,2.4]{Henk06PolRep}), which, together with the lower bound in Theorem \ref{main:thm}, determines
$\overline{N}(m,k)$ already up to a factor of $2$.

Theorem \ref{main:cor} will follow from Theorem \ref{main:thm}. The proof of Theorem~\ref{main:thm} consists of a geometric and a combinatorial part. In the geometric part, given in Section~\ref{geometric arguments}, we  show  that essentially it suffices to
consider unbounded polygons whose two unbounded edges are not parallel and express the minimal $n$ from Theorem~\ref{main:thm} for such polygons in purely combinatorial terms. In the combinatorial part, given in Section~\ref{combinatorial arguments}, we estimate this minimal $n$ by employing results on Gray codes (cf. \cite{MR1491049} and \cite{Knuth05}). 

By Theorem~\ref{main:cor}, $N(m,k) \sim \overline{N}(m,k) \sim m/k$ for $k \le m/\log_2 m$ as $m \rightarrow \infty.$ Let us give a few remarks for other choices of $k.$ We only discuss $N(m,k)$ here, since the comments on $\overline{N}(m,k)$ are analogous. Using results from real algebraic geometry, one can show that there exists a function $k_0(m)$ such that $N(m,k) = 2$ for every $k \ge k_0(m)$. This follows directly from the Theorem of Br\"ocker and Scheiderer \RAGbook[\S\,6.5,\,\S\,10.4] and  \cite[Theorem~5.1]{PrestelBooklet}, which are formulated in the framework of real closed fields. One should point out that the existence of $k_0(m)$ is a non-constructive result. In fact, on the one hand it is known that every $m$-gon $P$ can be described by two polynomials $p_1, p_2$, see \cite{Bernig98}. On the other hand, known formulas for $p_1, p_2$  depend on the metric structure of the $m$-gon $P$ so that the degrees of $p_1, p_2$ are not bounded from above for any fixed $m.$ Already for the case $m=5$ it is not clear how to (efficiently) construct polynomials $p_1, p_2$ representing an arbitrary given pentagon $P$ and having degrees bounded from above by an absolute constant, for further details see \cite{AveHenkDCG}. Finally, it would be interesting to estimate $N(m,k)$ for $m/\log_2 m \le k \le k_0(m)$, where $N(m,k)$ changes from about $\log_2 m$ to $2.$ 

As for $\mathcal{P}$-representations of polyhedra of arbitrary dimensions, we point out that recently it was shown \cite{AveBroe10}
that every $d$-dimensional polyhedron can be represented by $d$ non-strict polynomial inequalities. This provides a positive answer to a conjecture posed in \cite{MR2166533}. We also refer to
\cite{vomHofe}, \cite{Bernig98}, \cite{MR2166533},  \cite{AveHenkDCG}, \cite{AveHenkRepSimplePolytopes} for preliminary results,
 \cite{AveSomePropSemiAlg}, \cite{AveRepElemSemiAlg} for related results, and \cite{Henk06PolRep} for a survey on that topic.  For the case $d\ge 3$ currently no results on the interplay between the degrees and the number of polynomials in $\mathcal{P}$-representations of $d$-dimensional polyhedra seem available. 

\section{Geometric arguments} \label{geometric arguments}

Let $P$ be a polygon in $\real^2.$ By $\edges$ we denote the set of all edges of $P$. For every edge $I \in \edges$,
we fix a real polynomial $p_I(x)$ of degree one, vanishing on $I$ and non-negative on $P$.  

The following proposition is easily derived by standard arguments from real algebraic geometry.
For completeness we give a proof (see also \cite[Prop. 2.1]{MR1976602} and \cite{AveSomePropSemiAlg}).

\begin{proposition} \label{08.07.22,14:01} 
	Let $P$ be a polygon and let $p_1,\ldots,p_n$ be real polynomials such that $P=\{p_1\ge 0,\dots,p_n\ge 0\}$ or $\intr P = \rep{p_1}{p_n}{>}$. Then for every
$I \in \edges$ the polynomial $p_I$ is a factor of odd multiplicity of some  $p_i$, $1\le i\le n$.
\end{proposition} 

\begin{proof}[{\bf Proof}]
Fix an edge $I\in E$ and let $\prod_{i=1}^n p_i=p_I^{e}\prod_{j=1}^m q_j$ with $e\in\natur\cup\{0\}$ and irreducible polynomials
$q_j$ not divisible by $p_I$, $1\le j\le m$. Since every $q_j$ has only finitely many zeros in $\{p_I=0\}$ there exists an
$x\in I$ with $q_j(x)\not=0$ for all $j=1,\dots,m$, and hence, by continuity, there is an neighboorhood $\mathcal{N}$
of $x$ on which for all $j=1,\dots,m$ the sign of $q_j$ is constant.
If $p_I$ would have even multiplicity in every $p_i$ then the sign of every $p_i$ would be constant on
$\mathcal{N}\setminus\{p_I=0\}$ so that, in both cases $P=\{p_1\ge 0,\dots,p_n\ge 0\}$ and $\intr P = \rep{p_1}{p_n}{>}$,
either $\mathcal{N}\setminus\{p_I=0\}\subseteq P$ or $(\mathcal{N}\setminus\{p_I=0\})\cap \intr P=\emptyset$.
But this is impossible since $x$ is a boundary point of $P$ and
hence $\mathcal N \cap\{p_I<0\}\not\subseteq P$ and $\mathcal N \cap\{p_I>0\}\cap \intr P\not=\emptyset$.
\end{proof}

Recall that a real polynomial is called {\em squarefree} if every irreducible factor of it occurs with multiplicity one.
\begin{lemma}\label{lemma:squarefree}
Let $m,k\in\natur$. The quantity $N(m,k)$ (resp. $\overline{N}(m,k)$) from Theorem \ref{main:thm} is attained at
squarefree polynomials $p_1,\dots,p_n$ satisfying $(A)$ and $(C)$ (resp. $\overline{A}$ and $(C)$).
\end{lemma}
\begin{proof}[{\bf Proof}]
Given a polygon $P$ with $m$ edges and polynomials $p_1,\ldots,p_n$ satisfying the conditions ($A$) and ($C$), we can clearly
replace factors of $p_1,\ldots,p_n$ with odd multiplicities by the same factors with multiplicity one while maintaining the conditions (A) and (C).
Furthermore, factors
with even multiplicities can just be dropped out. In fact, if for some edge $I$, $p_I$ is a factor of $p_1$ with
even multiplicity, say with multiplicity two, then, in view of Proposition~\ref{08.07.22,14:01}, $p_I$ is a factor of odd multiplicity of some $p_i$, $1<i\le n$, say $p_2$, so that $\{p_1>0,p_2>0\}=\{p_1 / p_I^{2}>0,p_2>0\}$.

Similarly, given polynomials $p_1,\ldots,p_n$ satisfying the conditions ($\overline{A}$) and ($C$), odd multiplicities
can be replaced by multiplicity one, and factors with even multiplicities can be dropped out, since, with the notation above,
$P\subseteq \{p_1/p^2_I\ge 0, p_2\ge0\}\subseteq\{p_1\ge 0,p_2\ge 0\}$.

Performing the above reductions we obtain polynomials $p_1,\dots,p_n$ each of which being a product of at most $k$ different
$p_I$, $I\in E$. 
\end{proof}
For every set $F \subseteq \edges$ of edges of the polygon $P$ we put 
$$p_F(x) := \prod_{I \in F} p_I(x).$$

For $x \in \real^2 \setminus P$ we define the set $\edges_{<}(x) := \setcond{I \in \edges}{p_I(x) <0}.$ Geometrically, $\edges_{<}(x)$ is the set of edges of $P$ illuminated from $x.$ The sets $\edges_{\le }(x)$ and $\edges_{=}(x)$ are defined analogously.
\begin{lemma} \label{illum:crit}
	Let $P$ be a polygon and let $F_1,\ldots,F_n \subseteq \edges$. Then the following statements hold.
	\begin{enumerate}[I.]
		\item \label{08.07.22,13:59} If $P = \rep{p_{F_1}}{p_{F_n}}{\ge}$ then $\intr P = \rep{p_{F_1}}{p_{F_n}}{>}.$
		\item \label{08.07.22,14:00} If $\intr P = \rep{p_{F_1}}{p_{F_n}}{>}$ then $\edges = F_1 \cup \ldots \cup F_n.$ 
		\item \label{illum:crit:nonneg} $P=\rep{p_{F_1}}{p_{F_n}}{\ge}$ if and only if for every $x \in \real^2 \setminus P$ there exists $i \in \{1,\ldots,n\}$ such that $\card (F_i \cap \edges_{<}(x))$ is odd and $F_i \cap \edges_{=}(x) = \emptyset.$ 
		\item \label{illum:crit:pos} $\intr P = \rep{p_{F_1}}{p_{F_n}}{>}$ if and only if for every $x \in \real^2 \setminus P$ there exists $i \in \{1,\ldots,n\}$ such that $\card (F_i \cap \edges_{<}(x) )$ is odd.
	\end{enumerate} 
\end{lemma} 
\begin{proof}[\bf{Proof}]
The proof is straightforward. Part~\ref{08.07.22,14:00} follows directly from Proposition~\ref{08.07.22,14:01}. Let us show only
\ref{illum:crit:pos}. We start with the ``only if'' part.
Let $x \in \real^2 \setminus P.$ If $p_I(x) \ne 0$ for every $I \in \edges$ then, by assumption,
there exists $i \in \{1,\ldots,n\}$ such that $p_{F_i}(x) <0$ and consequently $\card(F_i \cap \edges_{<}(x) )$ is odd.
If $p_I(x)=0$ for some $I \in \edges$ then we can choose an $y \in \real^d \setminus P$
with $\edges_{<}(y)=\edges_{<}(x)$ and $p_J(y) \ne 0$ for every $J \in \edges$. Thus, by the previous case, there exists
$i \in \{1,\ldots,n\}$ such that $\card(F_i \cap \edges_{<}(x) )$ is odd. Next we show the ``if'' part. Since obviously
$\intr P\subseteq \{p_{F_1}>0,\dots,p_{F_n}>0\}$, we have to show that for every $x\in\real^2\setminus\intr P$ there exists
$i\in\{1,\dots,n\}$ with $p_{F_i}(x)\le 0$. If $x\notin P$ this is clear by assumption. If
$x\in P\setminus\intr P$ we have $x\in I$ for some $I\in E$, and we can choose an $y\in \real^2\setminus P$ with
$E_{<}(y)=\{I\}$. Then, by assumption, there exists $i\in\{1,\dots,n\}$ with $I\in F_i$, so $p_{F_i}(x)=0$.  
\end{proof} 

We remark that the converse of Lemma~\ref{illum:crit}.\ref{08.07.22,13:59} is not true in general. For example, the interior of the quadrant $P:=\setcond{x \in \real^2}{x_1 \ge 0, x_2 \ge 0}$ is equal to $$\setcond{x \in \real^2}{x_1 > 0, \ x_2 > 0}= \setcond{x \in \real^2}{x_1 > 0, \ x_1 x_2 >0},$$ and relaxing the inequalities on the right hand side gives
$\setcondbegin{x \in \real^2}{x_1 \ge 0,} \setcondend{x_1 x_2 \ge 0}$, the union of $P$ with the $x_1$-axis. 

By $n(P,k)$ (resp. $\overline{n}(P,k)$) we denote the minimal $n \in \natur$ such that there exist $F_1,\ldots,F_n \subseteq \edges$, each $F_i$ of cardinality at most $k$, such that $\intr P=\rep{p_{F_1}}{p_{F_n}}{>}$ (resp. $P=\rep{p_{F_1}}{p_{F_n}}{\ge}$).
Clearly, in view of Lemma~\ref{lemma:squarefree}, the quantity $N(m,k)$ (resp. $\overline{N}(m,k)$) in Theorem~\ref{main:thm} is the maximum of $n(P,k)$ (resp. $\overline{n}(P,k)$) over all polygons $P$ with $m$ edges.

We can now state our two main propositions in this section.
\begin{proposition} \label{08.07.22,15:12} Let $P$ be an unbounded polygon with $m$ edges such that the unbounded edges of $P$ are not parallel. Then the quantities $n(P,k)$ and $\overline{n}(P,k)$ depend only on $m$ and $k.$ Furthermore, $n(m,k):=n(P,k)$ and $\overline{n}(m,k):=\overline{n}(P,k)$ are determined as follows:
\begin{enumerate}[I.]
	\item \label{09.02.19,17:05} $n(m,k)$ is the minimal $n \in \natur$ such that there exist sets $S_1,\ldots, S_n \subseteq \{1,\ldots,m\}$ satisfying the following conditions:
	\begin{enumerate}
		\item[(I)] \label{n_1 odd cov} for every $a, b \in \natur$ with $1 \le a \le b \le m$ there exists $i \in \{1,\ldots,n\}$ such that $\card (S_i \cap \{a,\ldots,b\})$ is odd;
		\item[(K)]  \label{n_1 card bound} $\card S_i \le k$ for all $i \in \{1,\ldots,n\}$.
	\end{enumerate}
	\item \label{08.07.31,13:46} $\overline{n}(m,k)$ is the minimal $n \in \natur$ such that there exist sets $S_1,\ldots,S_n \subseteq \{1,\ldots,m\}$ satisfying the following conditions:
	\begin{enumerate}
		\item[(J)] \label{n_2 odd cov} for every $a, b \in \natur$ with $1 \le a \le b \le m$ there exists $i \in \{1,\ldots,n\}$ such that $\card (S_i \cap \{a,\ldots,b\})$ is odd and  $\{a-1,b+1 \} \cap S_i = \emptyset$;
		\item[(K)] \label{n_2 card bound} $\card S_i \le k$ for all $i \in \{1,\ldots,n\}$.
	\end{enumerate}
\end{enumerate} 
\end{proposition}
\begin{proof}[{\bf Proof}]
Let $I_1, I_2,\dots, I_m$ denote the edges of $P$ in consecutive order. Thus, $I_1$ and $I_m$ are the non-parallel unbounded
edges. By identifying edge sets $F_1,\dots, F_n\subseteq E=\{I_1,\dots,I_m\}$ with subsets
$S_1,\dots,S_n\subseteq\{1,\dots,m\}$ in the obvious way, part \ref{09.02.19,17:05} follows from Lemma~\ref{illum:crit}.\ref{illum:crit:pos} and
the equality
\[
	\setcond{\edges_<(x)}{x \in \real^2 \setminus P} = \setcond{\{I_a,\ldots,I_b\}}{1 \le a \le b \le m},
\]
and part \ref{08.07.31,13:46} follows from Lemma~\ref{illum:crit}.\ref{illum:crit:nonneg} and the equality
\begin{multline*}
	\setcond{\big(\edges_<(x), \edges_=(x)\big)}{x \in \real^2 \setminus P} = \\
 \setcond{\big(\{I_a,\ldots,I_b\},F\big)}{1 \le a \le b \le m,\, F\subseteq\{I_{a-1},I_{b+1}\}\cap E}.
\end{multline*}
\end{proof}
\begin{proposition} \label{08.11.12,12:28} For $m, k \in \natur$ let $N(m,k)$ and $\overline{N}(m,k)$ be defined as in Theorem~\ref{main:thm}, and $n(m,k)$ and $\overline{n}(m,k)$ as in Proposition~\ref{08.07.22,15:12}. Then one has
\begin{align*}
	n(m,k) & \le  N(m,k) \le 1+n(m-1,k), \\
	\overline{n}(m,k)  & \le  \overline{N}(m,k) \le 1+\overline{n}(m-1,k).
\end{align*}
\end{proposition}
\begin{proof}[{\bf Proof}]
	The lower bounds for $N(m,k)$ and $\overline{N}(m,k)$ are clear by Proposition \ref{08.07.22,15:12} and the remark
preceding it. As for the upper bounds consider an arbitrary polygon $P$ with $m+1$ consecutive edges $I_0,\ldots,I_m$.
Notice that 
	\begin{equation*}
		\setcond{E_<(x)}{x \in \real^2 \setminus P, \ I_0 \not\in E_<(x) } = \setcond{\{I_a,\ldots,I_b\}}{1 \le a \le b \le m}.
	\end{equation*}
	Consider sets $S_1,\ldots,S_n \subseteq \{1,\ldots,m\}$ satisfying {\itshape(I)} and {\itshape(K)}. As in the proof of
Proposition~\ref{08.07.22,15:12} we identify each set $S_i$ with the corresponding edge set $F_i \subseteq \{I_1,\ldots,I_m\}$.
In view of Lemma~\ref{illum:crit}.\ref{illum:crit:pos} we have $\intr P = \{p_{I_0} > 0, p_{F_1}>0,\ldots,p_{F_n}>0\}$.
Hence $N(m+1,k) \le 1+ n(m,k).$ The upper bound for $\overline{N}(m,k)$ is proved analogously and employs Lemma~\ref{illum:crit}.\ref{illum:crit:nonneg}. 
\end{proof}

\section{Combinatorial arguments and conclusion} \label{combinatorial arguments}

In view of the results from Section~\ref{geometric arguments}, the proofs of  Theorem~\ref{main:cor} and \ref{main:thm}
require the determination of the asymptotic behavior of the functions $n(m,k)$ and $\overline{n}(m,k)$
(defined in Proposition \ref{08.07.22,15:12}).
We will show the following. 
\begin{theorem}~\label{comb_prop}
Let $m, k \in \natur$, $k \le m.$ Then 
\begin{align*}
\max\{m/k, \log_2 (m\!+\!1)\} \le n(m,k) \le \overline{n}(m,k)\le
 \max \left\{\frac{m}{k} + 3 \log_2 \frac{m}{k} +C,\,\ceil{\log_2(m\!+\!1)}\right\},
\end{align*}
where $C>0$ is an absolute constant.
\end{theorem}

In particular,  for every $\eps>0$ and all $(1+\eps) m / \log_2 m \le k \le m$ with $m$ sufficiently large, we have the exact
values $\overline{n}(m,k)=n(m,k) = \ceil{\log_2(m+1)}$.

The proof of Theorem~\ref{comb_prop} employs results on Gray codes and is based on the following
reformulations for $n(m,k)$ and $\overline{n}(m,k)$.

For subsets $S_1,\dots,S_n$ of $\{1,\dots,m\}$ let $M$ denote their incidence matrix, i.e., the $n\times m$-matrix whose 
$(i,j)$-entry is $1$ if $j\in S_i$ and $0$ otherwise. Let $c_j$ be the $j$-th column of $M$, and let $M'$ denote the binary
$n\times(m+1)$ matrix whose $j$-th column is $c'_j$ is $c_1+\dots+c_j \modulo{2}$, $j=0,\dots,m$. In particular, $c'_0$ is the zero column. Clearly, $c_j=c'_{j-1}+c'_j \modulo{2}$ for all $j=1,\dots,m$.
Now it is easy to see that the conditions  $(I)$, $(J)$ and $(K)$ (in Proposition \ref{08.07.22,15:12})  on the
sets $S_1,\dots,S_n$ are equivalent to the following conditions $(I')$, $(J')$ and $(K')$, respectively, on the matrix $M'$:
	\begin{enumerate}
		\item[$(I')$] \label{gray cond 1} the $m+1$ columns of $M'$ are pairwise distinct;
	        \item[$(J')$] \label{gray cond 2} for every $a, b \in \natur$ with $0 \le a < b \le m$ there exists
                                                $i \in \{1,\ldots,n\}$ such that $M'_{i,a-1}=M'_{i,a}\neq M'_{i,b}=M'_{i,b+1}$
                                                (with $M'_{i,-1}:=M'_{i,0}$ and $M'_{i,m+1}:=M'_{i,m}$);
	        \item[$(K')$] \label{gray cond 3} every row of $M'$ has at most $k$ bit changes.
	\end{enumerate}
Here a {\em bit change in row $i$} of $M'$ is a column index $j\in\{1,\dots,m\}$ for which $M'_{i,j-1}\not=M'_{i,j}$.

Thus, $n(m,k)$ (resp. $\overline{n}(m,k)$) is the minimal $n\in\natur$ such that there exits a binary $n\times(m+1)$-matrix $M'$ starting with the zero column and satisfying $(I')$ and $(K')$ (resp. $(J')$ and $(K')$). 

Next we cite some facts about Gray codes (cf. \cite{MR1491049} and \cite{Knuth05}). Recall that an \emph{$n$-bit Gray code} is a binary
$n\times 2^n$-matrix $G$ with distinct columns $g_j$, $j=0,\dots 2^n-1$, such that $g_j$ and $g_{j+1}$ differ in exactly
one coordinate, $j=0,\dots,2^n-1$, where $g_{2^n}:=g_0$.
Without loss of generality we shall always assume that $g_0$ is the zero column.

A {\em bit run} of length $\ell$ of $G$ is a $1\times\ell$-submatrix $(G_{i,j+1},\dots,G_{i,j+\ell})$
of $G$ with $G_{i,j}\not=G_{i,j+1}=G_{i,j+2}=\dots=G_{i,j+\ell}\not=G_{i,j+\ell+1}$, where the columns are indexed modulo $2^n$.

The upper bound in Theorem \ref{comb_prop} essentialy follows from the following result due to Goddyn and Gvozdjak \cite{GoddynGvozdjak03}, see also
\cite[Section 7.2.1.1]{Knuth05}.
\begin{theorem} \label{gray code results} For every $n \in \natur$ there exists an $n$-bit Gray code for which every bit run has length at least $n-3\log_2 n$. \eop
\end{theorem}

In fact, as shown in \cite{GoddynGvozdjak03}, Theorem~\ref{gray code results} holds with the slightly stronger bound $\floor{n-2.001\log_2 n}$.

We will also need the following simple observation.
\begin{lemma}\label{J'lemma}
Every $n$-bit Gray code $G$ satisfies condition (J') (with $M'=G$).
\end{lemma}

\begin{proof}[{\bf Proof}]
Fix $a,b\in\mathbb{N}$ with $0<a<b<m$. Consider the row indices $j,k\in\{1,\dots,n\}$ for which $M'_{j,a-1}\not=M'_{j,a}$ and
$M'_{k,b}\not=M'_{k,b+1}$. Then, since $M'$ is a Gray code, $M'_{i,a-1}=M'_{i,a}$ for all $i\not=j$, and  $M'_{i,b}=M'_{i,b+1}$
for all $i\not=k$. In particular, the submatrix of $M'$ corresponding to the row indices $j,k$ and column
indices $a-1,a,b,b+1$ contains two equal columns (in both cases $j\not=k$ and $j=k$). Since the four columns of
$M'$ corresponding to the column indices $a-1,a,b,b+1$ are pairwise different, there exists an $i\not=j,k$ for which
$M'_{i,a-1}=M'_{i,a}\neq M'_{i,b}=M'_{i,b+1}$. 

A similar argument works if $a=0$ or $b=m$.
\end{proof}

\begin{proof}[\bf{Proof of Theorem~\ref{comb_prop}}] We shall use the matrices $M$ and $M'$ introduced above.

The lower bound on $n(m,k)$ is clear in view of the conditions on $M$ and $M'$. Indeed, with $a=b$ in
condition $(I)$ we see that every element from $\{1,\dots,m\}$ is covered by at least one  of the sets $S_1,\dots,S_n$,
which in view of $(K)$ yields $n\ge m/k$. (Alternatively, we can count the total number of bit changes in $M'$.)
Further, we have $n\ge\log_2(m+1)$ since the $m+1$ columns of $M'$ are distinct.

The bound $n(m,k) \le \overline{n}(m,k)$ is trivial.

It remains to show the upper bound on $\overline{n}(m,k)$. Fix the minimal $n \in \natur$ satisfying $m\le 2^n-1$ and $m \le k(n-3\log_2n)$ 
and consider an $n$-bit Gray code $G$ as in Theorem~\ref{gray code results}. Then the
$n\times(m+1)$-matrix $M'$ consisting of the first $m+1$ columns of $G$ starts with the zero column and clearly satisfies
condition $(K')$. By Lemma \ref{J'lemma}, $M'$ also satisfies condition $(J')$. Thus, $n(m,k)\le n$, from which the desired upper
bound follows by elementary calculations. Indeed, if $n>\ceil{\log_2(m+1)}$ then, by minimality of $n$,
$m/k>((n-1)-3\log_2(n-1))\ge (n-1)/4$ (the last inequality for $n\ge 8$) and thus
\[
n-1<m/k+3\log_2(n-1)<m/k+3\log_2(4m/k)=m/k+3\log_2(m/k)+6.
\]
\end{proof}

\begin{proof}[\bf{Proof of Theorem~\ref{main:thm}}]
The inequality $N(m,k)\le\overline{N}(m,k)$ is clear in view of Proposition \ref{illum:crit}.\ref{08.07.22,13:59}.
The remaining inequalities follow from Proposition~\ref{08.11.12,12:28} and Theorem~\ref{comb_prop}. 
\end{proof}

\begin{proof}[\bf{Proof of Theorem~\ref{main:cor}}] The lower bounds on $N(m,k)$ and $\overline{N}(m,k)$ follow from
Proposition~\ref{08.07.22,14:01}. The upper bounds follow directly from Theorem~\ref{main:thm}.
\end{proof}

\section*{Acknowledgement}
We thank Ludwig Br\"ocker, Martin Henk, Bettina Matzke and Claus Scheiderer for helpful discussions.
\small

\bibliographystyle{amsalpha}

\bigskip
 \begin{tabular}{l}
        \textsc{Gennadiy Averkov, Christian Bey,
	} \\ \textsc{Universit\"atsplatz 2,} \textsc{Fakult\"at f\"ur  Mathematik,} \\
        \textsc{Otto-von-Guericke-Universit\"at Magdeburg,} \\ \textsc{D-39106 Magdeburg}	\\
        \emph{e-mails}: \texttt{\{averkov,bey\}@ovgu.de}\\
	\emph{web}: \texttt{http://fma2.math.uni-magdeburg.de/\{$\sim$averkov,$\sim$bey\}} \\
    \end{tabular} 

 \end{document}